\input allknots.sty                                      \advance\vsize by10pt

\^    Spacefilling knots. (November 9, 1994; final version: November 27,1996)

\*    _Peter Schmitt_


\vfill

0. Summary.
! This note describes how to construct toroidal polyhedra
  which are homotopic to a given type of knot
  and which admit an isohedral tiling of 3-space.

* 1. Introduction.

A _knotted_ (toroidal) polyhedron (of a certain type) is a polyhedron
  which is homeomorphic to a closed torus
  and homotopic to a knot of this type.
This note describes in detail the construction of knotted polyhedra
  which are spacefillers, i.e., which admit a tiling of 3-space.
  The type of the knot may be chosen arbitrarily,
  and the resulting tiling is both periodic and isohedral,
  i.e., the group of symmetries acts transitively on the tiles.
This construction is a generalization of that used in~[4] and~[5].
  It has been briefly indicated in~[6].

These investigations were initiated by a lecture of Schulte~[7] during
  which he suggested to search for a spacefilling knotted polyhedron.
  (However, the published version~[8] does not mention knotted tiles.)
Other constructions of spacefilling knots
  (independently, but also initiated by Schulte)
  are given by Kuperberg~[3] and Adams~[1].
  Both authors describe knotted polyhedra which admit
       a lattice tiling of 3-space, without interlocking.
  Kuperberg's approach is simplified by Zaks~[9].

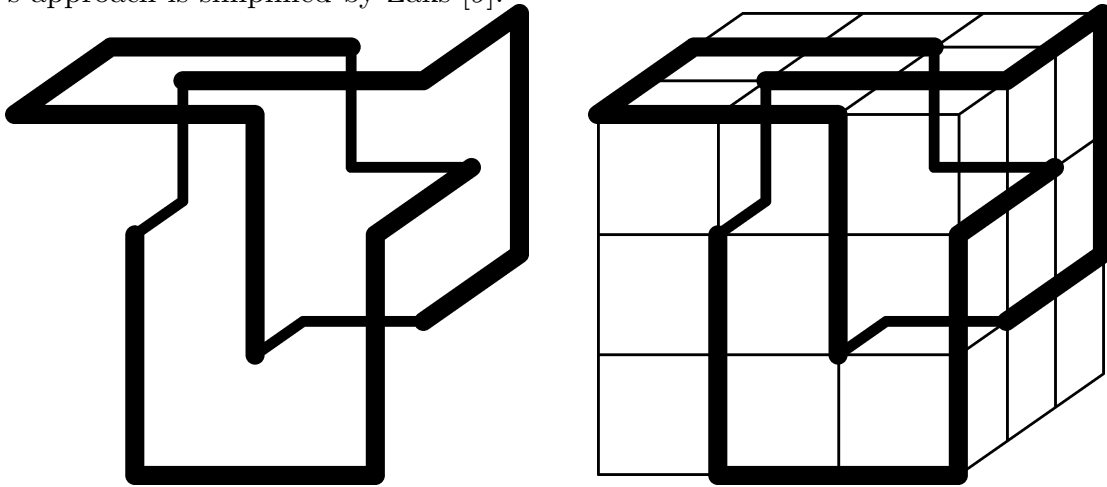
\Fig 1: A trefoil knot and its embedding into a cubic grid.

* 2. An outline of the construction.

Consider the graph in 3-space formed by the edges of a cube lattice.
Obviously, it is possible to embed a knot of any desired type into this graph
  (Fig.~1 shows one of the possible realizations of the simplest case:
  a trefoil knot embedded into a $3\times3\times3$ cubic grid.)
Apply all the translational symmetries of the lattice to obtain a
  family of congruent copies of embedded knots. This family covers the graph.
  Duplicate the family by translating the vertices to the centres of the cubes.
Intuitively, it is clear that in three dimensions there is enough freedom to
  modify the embedded knot slightly to obtain a knot
  (homotopic and `close' to the original one,
  but no longer contained in the edge graph)
  such that its copies by translation are pairwise disjoint.
Finally, blow up the knot to a topological torus such that
  its congruent copies form a tiling of 3-space,
  while taking care that the homotopy class remains unchanged by this process.

The procedure just outlined can, of course, be realized in several ways.
  One specific construction is described in the next section.
  It is essentially straightforward,
  but the details are partially rather cumbersome to write down explicitly.

* 3. The construction in detail.

ø Step 1: A partition of 3-space

Consider the tiling \T of 3-space by
  unit cubes with centres at the elements of the integer lattice \L.
For the purpose of the construction this tiling is split into
  two families, $T1$ and $T2$, of tiles
  according to the majority of even or odd integers among the coordinates
  of the centres.
The unions $S1$ and $S2$ of these two families are congruent, have disjoint
  interiors, and (together) cover space.
(More formally,
  $  T := \{ C+p | p \in L  \} $,
  $ Ti := \{ C+p | p \in Li \} $ and
  $ Si := \bigcup Ti = C+Li $
    where $ C := [-1/2,1/2]^3 $, $ L = \Z^3 $, and
\\  $ L1 := \{ p \in L | `at most one coordinate odd' \} \subset L $ and
\\  $ L2 := \{ p \in L | `at most one coordinate even' \} = L1+(1,1,1) $.
\\  This partition of 3-space into two sets has also been used by Debrunner~[2].)

Therefore, if $S1$ can be tiled monohedrally by some prototile \K,
  then $\R^3$ can be tiled by it, as well.
Furthermore, note that $S1$ and $S2$ are homotopic to the
  3-dimensional grid \G (of lines parallel the axes)
    $$ G := \{ p \in \R^3 | `at least two coordinates are integers' \} $$
  generated by the integer lattice or, equivalently,
  to $ 2G \subset S1 $.
Guided by the geometric structure of \G,
  the set $L1$ can be divided into two classes,
  the _vertices_
     $$ V := \{ p \in L | `all coordinates even' \} = 2L \subset L1 $$
  and the _edges_ (represented by their midpoints)
     $$ E := \{ p \in L | `two coordinates even and one odd' \} \subset L1 $$
  (with three subclasses $Ex$, $Ey$, and $Ez$,
      according to the odd coordinate).

ø Step 2: Embedding of the knot

Every knot
  (or more precisely, every knot that has a piecewise linear representation)
  can be embedded into $ 2G \subset S1 $
  by a closed path \P consisting of points of $L1$
  and line segments (each parallel to an axis) between these points.
Note that, among the lattice points along such a path,
  vertices and edges necessarily alternate.

This path \P will be the basis for the construction of
  the desired knotted spacefiller \K, a polyhedral torus homotopic to \P.
Consider all congruent copies $ P+p $ of \P obtained
  by adding a point $ p \in V = 2L $.
Obviously, these translates cover \G, they are, however, far from disjoint.
It is therefore necessary to modify \P to obtain an (equivalent) knot \Q
(homotopic to \P, both in $\R^3$ and in $S1$) for which the translates are
pairwise disjoint.

The path (or knot) \P can be described by the sequence
   $ v0, e1, v1, e2, \dots en, vn = v0 $
   of its vertices and edges.
In order to facilitate the construction,
  set $i(x)=n$ and assume that
  $ en      = (-1,0,0) $,
  $ vn = v0 = ( 0,0,0) $,
  $ e1      = ( 1,0,0) $,
  and that there are indices $i(y)$ and $i(z)$
  such that both $e_{i(y)}$ and $e_{i(y)+1}$ belong to the class $Ey$,
        and both $e_{i(z)}$ and $e_{i(z)+1}$ belong to the class $Ez$.
For example, an admissable representation of the trivial knot would be
  (omitting the edges):
  $$\eqalign{
  v_0 = (0&,0,0)-(2,0,0)-(2,2,0)-(2,4,0)-(2,4,2)-(2,4,4)-(0,4,4)-        \cr
         & -(-2,4,4)-(-2,2,4)-(-2,0,4)-(-2,0,2)-(-2,0.0)-(0,0,0) = v_{12}
  \cr}$$

ø Step 3: Modification of the knot

Since the translations by elements of \V act transitively on the vertices
   and on each of the classes $Ex$, $Ey$, $Ez$ of edges,
in order to construct \Q it is sufficient to consider the
   four cubes \C, $Cx$, $Cy$, and $Cz$
       (where $ Cx = C+(1,0,0) $, etc.):
   If the knots $ Q-vi $ ($ i=1,\dots,n $) do not intersect
   in these four cubes, then all the translates of \Q are pairwise disjoint.

To each vertex $vi$ assign a point $ pi = (x_i,y_i,z_i) \in C $
   such that (for $ i=1,\dots,n $) no coordinate occurs twice.
For technical reasons (whose purpose will become clear later)
   some additional assumptions on the $pi$ are convenient
   (and can obviously be satisfied):
   To $v0=vn$ assign $ p0=pn =p_{i(x)}= (0,0,0) $ and choose
   $x_{i(y)}$, $z_{i(y)}$, $x_{i(z)}$, and $y_{i(z)}$ such that
   their absolute values are less than that of all the other coordinates.
Furthermore, assume that $pi$ has been chosen such that
   (if applicable)
   the right angle determined by the two lines meeting in $pi$
   is intersected by the coordinate axis which is orthogonal to the plane
   determined by them.
(Applied to the sample path of Step~2 this means that, e.g.,
   for $p1$ the coordinate $x_1$ has to be chosen from the interval $(0,1/2)$
        and the coordinate $y_1$ has to be chosen from $(-1/2,0)$.)

When constructing \Q, the segment $ C \cap (P-vi) $
   --- which consists of two line segments (parallel to two of the axes)
       meeting in (0,0,0) ---
   is replaced by a corresponding segment $ C \cap (Q-vi) $.
      which is formed by two line segments
      (parallel to the original ones, and having the same directions)
      meeting in $pi$.
Therefore, for every point of $ C \cap (Q-vi) $
   at least two coordinates coincide with those of $pi$.
   By the choice of these coordinates (as pairwise distinct) this means
   that (for $ 1 \le i \le n $) these segments are pairwise disjoint.

For each edge $ ei \in Ex $
    either $ ei - v_{i-1} \in Cx $
        or $ ei - vi      \in Cx $
    (depending on the orientation).
To obtain the segment $ Cx \cap (Q-v_{i-1}) $, in the first case,
   the segment $ Cx \cap (P-v_{i-1}) $ has to be replaced by some path
   leading from $ (1/2,y_{i-1},z_{i-1}) $ to $ (3/2,y_i,z_i) $.
   A possible choice is the following path obtained by joining
   five points by four straight line segments
   $ (1/2  ,y_{i-1},z_{i-1}) $ --
   $ (1+x_i,y_{i-1},z_{i-1}) $ --
   $ (1+x_i,y_i    ,z_{i-1}) $ --
   $ (1+x_i,y_i    ,z_i    ) $ --
   $ (3/2  ,y_i    ,z_i    ) $.
In the second case,
   the replacement $ Cx \cap (Q-vi) $ of $ Cx \cap (P-vi) $
   is obtained similarly by joining the points
   $ (1/2  ,y_i    ,z_i    ) $ --
   $ (1+x_i,y_i    ,z_i    ) $ --
   $ (1+x_i,y_{i-1},z_i    ) $ --
   $ (1+x_i,y_{i-1},z_{i-1}) $ --
   $ (3/2  ,y_{i-1},z_{i-1}) $.
In most cases, these paths do not intersect
   by the choice of coordinates (just as in the case of vertices):
For each point of the segment corresponding to $ei$ (in $Cx$)
  (at least) two of the coordinates coincide
  with those of either $pi+(1,0,0)$ or $p_{i-1}+(1,0,0)$.
  Therefore, the segments corresponding to two edges are disjoint
  if they do not have a common endpoint (vertex).
However, if a vertex $vi$ occurs twice
  then both $ei$ and $e_{i+1}$ belong to $Ex$ and,
  if these edges are positively oriented,
  the last part of the segment corresponding to $ei$ and
  the first part of the segment corresponding to $e_{i+1}$
  (or vice-versa, if they are negatively oriented)
  have the same $x$-$y$-coordinates and may intersect:
The line segment
  $ (1+x_i,y_i,z_i) $ -- $ (3/2,y_i,z_i) $
        (from $Q-v_{i-1}$ corresponding to $e_i$) intersects the line segment
  $ (1/2,y_i,z_i) $ -- $ (1+x_{i+1},y_i,z_i) $
        (from $Q-v_i$ corresponding to $e_{i+1}$)
  if $ x_i < x_{i+1} $ (for positively oriented edges), and similarly,
  $ (1/2,y_i,z_i) $ -- $ (1+x_i,y_i,z_i) $
        (from $Q-v_i$ corresponding to $e_i$) intersects
  $ (1+x_{i+1},y_i,z_i) $ -- $ (3/2,y_i,z_i) $
        (from $Q-v_{i+1}$ corresponding to $e_{i+1}$)
  if $ x_i > x_{i+1} $ (for negatively oriented edges).
But when constructing the segments corresponding to the edges in $Ex$
  the actual values of the $x_i$ are not important --- it is only used
  that they are distinct.
  (In fact, they could be replaced by arbitrarily chosen values.)
Therefore, this exceptional case can easily be handled
  by exchanging $x_i$ and $x_{i+1}$.
  Similarly, when three or more consecutive edges belong to $Ex$
  then, for the corresponding segments, the $x_i$ can be permuted
  to produce (descending or ascending) order.

For a unified notation, write
   $ (1/2,ai,bi) $ --
   $ (ci,ai,bi) $ --
   $ (ci,Ai,bi) $ --
   $ (ci,Ai,Bi) $ --
   $ (3/2,Ai,Bi) $
   to denote the path segment
   $ Q-v(i) $ (where $ v(i) = vi $ or $ v(i) = v_{i-1} $)
   corresponding to $ ei \in Ex $.

(The edges in $Ey$ and $Ez$ are, of course, treated similarly.)

ø Step 4: Thickening of the knot

For the final step
  divide the cubes $C$, $Cx$, $Cy$, and $Cz$ into sufficiently small cubes,
  i.e., cubes with a side length less than half the minimal distance
  occurring between the coordinates $x_i$, $y_i$, $z_i$, $-1/2$, and $1/2$.
This guarantees that each small cube meets at most one of the
  knot segments $Q-vi$.

Now construct $K0-vi$ by taking the
  union of all small cubes that meet $Q-vi$
  The result is a polyhedral torus $K0$ which is homotopic to \Q and \P
  and whose translates (by the elements of \V) still are pairwise disjoint.
`Blowing up' (or thickening) $K0$ so that
  the resulting body and its translates tile $S1$
  will complete the construction and provide a knotted polyhedron \K as desired.
  This task can be achieved as follows:

ø Step 4a: Edge segments

Consider first $Cx$.
If $(1/2,ai,bi)$ belongs to a segment $Q-v(i)$ then add to $K0-v(i)$
   all cubes meeting the line segment
   from $(1/2,ai,bi)$ to $(1/2,ai,\pm1/2)$
        (the sign depending on the sign of $bi$).
   Then add all the remaining cubes meeting the face $x=1/2$ of $Ex$
   --- a `two-sided comb' extending from (and including) the line
      $ (1/2,-1/2,z_{i(x)}) - (1/2,1/2,z_{i(x)}) $ ---
   to $Q-v_{i(x)}$.
   This works since the coordinates $ai$
        (taken from the $y_i$) are all different,
   and the coordinates $bi$ (from the $z_i$), for $ i \not= i(x)$,
   all have absolute values greater than $z_{i(x)}$.
Proceed analogously with the opposite face ($x=3/2$).

These two face layers (and analogously constructed ones in $Cy$ and $Cz$)
   serve to hide the interior of $Cx$ (or $Cy$ and $Cz$, respectively).
   By their construction,
   and because of the assumptions on the coordinates of $pi$,
   the segments (of $K0$),
   corresponding to two edges of \Q which meet orthogonally in $vi$,
   are extended only outwards or sidewards, but not inwards,
   and therefore they will not touch when the polyhedron \K is put together.

Next consider the values $ci$ (for edges in $Ex$),
   i.e., those values of the $x$-coordinate
   where, with the help of two orthogonal line segments,
   a change from one line parallel to the $x$-axis to another one takes place.
   For each $ ej \in Ex $ ($ j \not= i $) the path segment of $Q-v(j)$ passes
      either through $(ci,aj,bj)$ or $(ci,Aj,Bj)$.
   Accordingly, add to $K0-v(j)$ all cubes meeting the line segment
      $ (ci,aj,bj) - (ci,aj,\pm1/2) $
      (or the line segment $ (ci,Aj,Bj) - (ci,Aj,\pm1/2) $, respectively)
        where the sign is chosen such that this line
        does not meet the two orthogonal line segments (of $Q-v(i)$ at $x=ci$).
   Add all remaining cubes meeting
        the square cross section determined (in $Ex$) by $x=ci$
        (a simply connected set) to $K0-v(i)$.

So far, from the small cubes of $Cx$
   two layers at two faces and some layers (less than $n$) between them
   have been assigned to segments $K0-v(i)$
   In the sections between these layers $Q-v(i)$ (if in $Cx$)
   is a line segment parallel to the $x$-axis
   (and therefore $K0-v(i)$ is a prism parallel to this axis).
   In each section, extend these finitely many prisms
     (less than $n$, but at least one because of $e_{i(x)}$)
   to layers (parallel to the $x$-$z$-plane and of appropriately chosen width)
   which fill that section.
   (Since all the coordinates $y_i$ are different, no conflict occurs.)

The resulting segments $K-v(i)$
   all are simply connected and form a tiling of $Cx$.
   Proceed analogously for $Cy$ (extension again parallel to the $z$-axis)
   and $Cz$ (extension parallel to the $y$-axis).

ø Step 4b: Vertex segments

Now consider \C and use a construction similar to the previous one:
  The segment of $Q-v0=Q$ (corresponding to $v0$) passing through \C
  coincides with the $x$-axis since $p0=(0,0,0)$.
   If a segment of $Q-vi$ intersects the $x$-$y$-plane in $(x_i,y_i,0)$
   then add to $K0-vi$ all cubes meeting the line segment
       from $(x_i,y_i,0)$ to $(x_i,\pm1/2,0)$
       where the sign of $1/2$ equals the sign of $y_i$.
   (Again, no conflict occurs, since the $x_i$ are all different.)
   Then add all remaining cubes meeting the square with $z=0$ to $K0-v0$.
In other words:
  Add to $K0$ prisms parallel to the $y$-axis (in both directions)
     either until a face of \C is reached, or
     until another segment $K0-vi$ is met.
  In the latter case, add cubes to this segment (in the same way),
  outwards and parallel to the $y$-axis, until the face of \C is reached.
The result of this procedure is a square layer (at $z=0$) of small cubes,
   where each cube is assigned to one of the segments $K0-vi$.
Now, similarly, add more cubes to the extended segments $K0-vi$:
   To every cube of the square layer just constructed
   add cubes to form prisms parallel to the $z$-axis
   either until a face is reached,
   or until a cube of another segment $K0-vj$ is met
   --- in the latter case add the remaining cubes to this segment.

Since any right angle formed by a segment contains a coordinate axis
   these angles are enlarged
   only outwards or sidewards, but never inwards
   --- therefore the construction does not create loops,
   and the result is a a tiling of \C
       by simply connected polyhedra $K-vi$.
   Moreover, for the same reason no loops will occur (at an edge cube)
   when, finally, the segments are glued together
   to form the polyhedron \K.

ø The result

In the end, each one of the small cubes in $C \cup Cx \cup Cy \cup Cz $
   is assigned to one of the segments $K0-vi$ to form a segment $K-vi$,
   therefore the $n$ segments $K-vi$ form a tiling of this set, and
$$ K = \bigcup_i ((K-vi)+vi) \supset Q $$
   is a knotted prototile with the desired properties:
   Because of its construction this polyhedron
   belongs to the same homotopy class as \Q (and \P),
   and the family of its translates (by the elements of \V)
   forms a tiling of $S1$.

* 4. Remarks.

(a) _(Step~1)_ The construction does not depend on \P being a knot.
    It may as well start from an 1-dimensional manifold of higher genus,
    and may even posess more than one connected component
    (e.g., it may be a link).
    In this way, by starting from
    \P and some (disjoint) congruent copies of \P,
    some additional symmetries for the resulting tiling may be obtained.

(b) _(Step~2)_ Both the choice of the embedded path \P
    and the type of modifications leading to \Q influence
    how the copies of \K are linked in the final tiling.

(c) _(Step~1)_ Instead of $S1$ and $S2$ the construction can also be based
    on the `linked double-layers' used in~[3],
    or some similar divisions of 3-space.

* 5. References.

[1] Colin C. Adams, _Tilings of space by knotted tiles._
    The Math. Intelligencer 17 (2), 41--51 (1995).

[2] Hans E.~Debrunner, _Tiling three-space with handlebodies._
    Studia Scient.\ Math.\ Hung.~21 (1986), 201--202.

[3] Wlodzimierz Kuperberg, _Knotted lattice-like space fillers._
    Discrete Comput. Geom. 13, 561--567 (1995).

[4] Peter Schmitt, _A spacefilling trefoil knot._
    Note, June~14, 1993
    (to appear: \"Osterreich. Akad. Wiss., math.-naturw. Kl., Anzeiger).

[5] Peter Schmitt, _Another space-filling trefoil knot._
    Discrete Comput. Geom. 13, 603--607 (1995).

[6] Peter Schmitt, _Spacefilling knots._ Preprint, 1994.

[7] Egon Schulte, _Spacefillers of higher genus._
    Lecture, Austrian conference on discrete geometry. May 2--8, 1993.
    Neuhofen/Ybbs, Austria.

[8] Egon Schulte, _Space fillers of higher genus._
    J. Comb. Theory A 68, 438--453 (1994).

[9] Joseph Zaks, _Monohedrally knotted tilings of the 3-space (a short note)._
    Proc.~Combinatorics '94 (to appear).
                                                                   \vskip1cm
                                               \bigskip\obeylines\parskip0pt

Peter Schmitt, _Institut f\"ur Mathematik, Universit\"at Wien_
Strudlhofgasse 4, A--1090 Wien, Austria
|Peter.Schmitt@univie.ac.at|

\bye